\documentclass[12pt,amsmaths,preprinnumbers]{article}
\usepackage{psfig}
\usepackage{graphicx}
\usepackage{epsfig}
\usepackage{dcolumn}
\usepackage{bm}
\begin{document}
\newcommand{\be}{\begin{equation}}
\newcommand{\ee}{\end{equation}}
\newcommand{\ben}{\begin{eqnarray}}
\newcommand{\een}{\end{eqnarray}}
\newcommand{\bF}{\begin{figure}}
\newcommand{\eF}{\end{figure}}
\title{Turbo Codes over the Real Field}
\author{Animesh Datta \footnote{E mail: animesh@unm.edu}\\
{\it\small Department of Electrical Engineering, Indian Institute
of Technology, Kanpur, 208 016}\footnote{Present Address: Deptt. of Physics and Astronomy, University of New Mexico, Albuquerque, New Mexico, 87131, USA}}
\maketitle
%\begin{document}
\begin{abstract}
{\it In this article we extend the idea of Turbo codes onto the Real Field. The channel is taken to result in block erasures and the only noise as being that due to quantization. The decoding in this case is reduced to reconstruction of the lost values. The encoding is done using critically sampled filter banks and introduction of an interleaver is found to reduce the mean square quantization error drastically. The permutation that gives the best recoverability is obtained in the 2 Channel case. Results are also obtained for M channel case. The algorithm for reconstruction of the lost values in the absence of quantization noise is obtained.}

\end{abstract}
\newpage

\section{Introduction}
Novel communication and information services are being intoduced almost daily and demands for higher data rates and communication capacity continue to grow. This spectacular progress of communication is to a great extent due to consistent performance increase and cost reductionof the devices and circuit technology. Such advancements have also been fostered by major theoretical developments. The synergy components and signal processing techniques is considered to be the main cornerstone of modern communication equipments.
	
Since Shannon established the fundamental limits on the transmission rates of a digital communication system and motivated the search for coding techniques that approach the capacity, various coding techniques have been proposed. Convolutional codes did better than block codes that formed the art of coding and cryptography. The landmark development in this regard has been by Berrou, Glavieux and Thitimajshima when they proposed turbo error control codeing over finite fields ($GF(2)$) by which the gap between the capacity limit and practically feasible channel utilisation is almost closed.

In this article we to extend the idea of turbo codes proposed in \cite{turbo1},\cite{turbo2} to codes over real fields. 

The plan of the article would be as follows. In Section 2, we introduce turbo codes over $GF(2)$. In Section 3 we familiarise the reader with the techniques of frame theory and the theory of filter banks as convolutional encoders. In Section 4 we present our results on the encoding of real data using our {\it turbo encoders}. We compare them with normal oversampled filter banks and note the advantages. In Section 5 we present our results about the optimal permutation and also deal with M channel Encoders. Finally we conclude in Section 6 with a summary and a statement of some open problems.  
 
\section{Turbo Codes over $GF(2)$}

Turbo codes exploit an idea of connecting two codes and seperating them by an interleaver. The primary achieve a low error rate with an overall decoding complexity lower than that of a single code of corresponding performance. The low complexity is attained by decoding each component seperately. As the inner decoder generates burst errors,an interleaver is typically incorporated between the two codes to decorrelate the received symbols affected by burst errors. Another application of concatenation is using a bandwidth efficient trellis code as an inner code or concatenating two convolutional codes. 

The difference between turbo and serially concatenated codes is that in turbo codes two identical systematic component codes are connected in parallel. The information bits for the second code are not transmitted thus increasing the code rate relative to a corresponding serial concatenated code. The primary reason for using a long interleaver in turbo codes is to generate a concatenated code with large block length which leads to a large coding gain. The decoder consists of two concatenated decoders of the component codes seperated by the same interleaver. The component decoders are based on {\it max a posteriori} (MAP) probability algorithm or a {\it Soft Output Viterbi Algorithm} (SOVA) generating a weighted soft estimate of the input sequence. The iterarive decoding process performs information exchange between the two component decoders. By increasing the number of iterations in turbo-decoding, a bit error rate (BER) as low as $10^{-5}$ to ${10^{-7}}$ can be achieved at Signal-to-noise (SNR) close to the Shannon Capacity limit. 

\subsection{Turbo Encoder}

A turbo encoder is formed by parallel concatenation of two {\it Recursive Systematic Convolutional} (RSC) encoders seperated by a random interleaver. Thus turbo codes are also regarded as {\it Parallely concatenated convolutional codes} (PCCC).

\subsection{Interleaver}
The interleaver in turbo codes is a pseudo-random block scrambler defined by a permutation of $N$ elements with no repetition. The first role of the interleaver is to generate a long block code from small memory convolutional codes. Secondly, it decorrelates the inputs the two decoders so that an iterartive suboptimum decoding algorithm based on information exchange between the two component decoders can be applied. If the input sequence to the two components of the decoder are decorrelated, there is high probability that after the correction of some of the errors in one decoder some of the remainging errors should be correctable in the second decoder. In a pseudo random interleaver a block of $N$ inputs is read into the interleaver and read out pseudo randomly. The pseudo random interleaving pattern \emph{must} be available at the decoder as well.

\subsection{Decoder: Iterative Decoding}  

Turbo codes can be decoded by MAP or ML decoding methods based on the overall code trellis. These decoders can be implemented for small interleavers only as they are too complex for medium and large interleaver sizes. The practical importance of tubo codes lies in the availability of a simple suboptimum decoding algorithm. 

The iterative turbo decoding consists of two component decoders serially concatenated via an interleaver, identical to the one in the figure. The first MAP decoder takes as input the received information sequence ${\bf r_0}$ and the received parity sequence generated by the first encoder ${\bf r_1}$. The decoder then produces a soft output , which is interleaved and used to produce an improved estimate of the a priori probabilities of the information sequence if the second encoder.

The other two inputs to the second MAP decoder are the information sequence ${\bf \tilde{r}_0}$ and the received parity sequence produced by the second encoder ${\bf r_2}$. The second MAP decoder also produces a soft output which is used to improve the estimate of the a priori probabilities of the information sequence at the input of the first MAP decoder. The decoder performance can be improved by the iterative operation relative to a single operation serial concatenated decoder. The feedback loop is a distinguishing feature of the decoder and the name turbo code is given with reference to the turbo engine.

After a certain number of iterations the soft outputs of both MAP decoders stop producing further perfomance improvements. Then the last stage of the decoding makes the hard decision after deinterleaving.     

\section{Frames and Filter Banks}

\subsection{Frames}
In this chapter, we define discrete frames and summarize a few of their properties. For a detailed treatment of frames and generalized frames, the reader is referred to \cite{daub},\cite{kaiser}.

We consider frames in $K$ dimensional real space, i.e ${R}^{K}$. A set of $K$-dimensional vectors $\Phi_{\it F}\equiv{\varphi_{k}}, k=1,\dots,N$ is called a frame is there exist $A>0$ and $B < \infty$ such that 
\be
A\parallel {\bf z} \parallel ^2 \leq \sum_{k=1}^{N}|\langle{\bf z},\varphi_k\rangle|^2 \leq B \parallel {\bf z} \parallel ^2 
\ee  \label{frame}
for all $ {\bf z}\in {R}^{K}$ 

where $|\langle{\bf z},\varphi_k\rangle|$ denotes the inner product of ${\bf z}$ and $\varphi_k$, and $\parallel {\bf z} \parallel $ denotes the Euclidean norm of ${\bf z}$. $A$ and $B$ are called the frame bounds. The inner product $|\langle{\bf z},\varphi_k\rangle|$ gives the $k$th frame expansion coefficient of {\bf z}. Any finite set if vectors that spans ${R}^K$ is a frame. Therefore a frame will always have $ N\geq K$. The ratio $N/K$ is normally referred to as the redundancy of the frame.

A frame is called \emph{uniform} if each frame vector has unity magnitude,i.e., $\parallel{\varphi}_k\parallel=1$ for $k=1,2,\dots,N$. The frame $\Phi_F$ is associated with a frame operator $F$ which is defined as follows:
\be
F{\bf z}\equiv |\langle{\bf z},\varphi_k\rangle|\ee for all  $k=1,2,\dots,N$.

Therefore the frame expansion coefficients of {\bf z} are given by $F{\bf z}$. Using the frame operator $F$,the frame condition in eqn.\ref{frame} can be rewritten as 
\be
A\parallel {\bf z} \parallel ^2 \leq {\bf z}^tF^tF{\bf z} \leq B \parallel {\bf z} \parallel ^2
\ee for all ${\bf z}\in {R}^{K}$ where ${\bf z}^t$and $F^t$ denote the matrix transposes of ${\bf z}$ and $F$. It can be easily shown that the sum of the eigenvalues of $F^tF$ is equal to the sum of the squared lengths of the frame vectors.

A frame is called \emph{tight} if its bounds are equal,i.e., $A=B$.Therefore $\Phi_F$ is tight if and only if $F^tF=AI_k$. This implies that the columns of $F$ are orthogonal. In this case, $F^tF$ has eigenvalue A with algebraic multiplicity $K$. A frame is called \emph{snug} if $A\approx B$.	

With $N>K$, $F{\bf z}$ provides an over-complete expansion of the signal vector {\bf z}. In this case, the frame vectors cannot be orthogonal and therefore the dual frame of $\Phi_F$ is required to reconstruct the vector {\bf z} from its frame expansion coefficients. The dual frame \footnote{It is also referred to as the reciprocal frame by some authors.}of $\Phi_F$ is another frame defined as $\Phi_{\tilde F}\equiv {\tilde\varphi_k}$,$k=1,\dots,N$, where $\tilde\varphi_k=(F^tF)^{-1}\varphi_k$, for $k=1,2,\dots,N$. Since $F^tF$ is positive definite, it is invertible and hence the above definition is valid. It is easily seen that the frame operator associated with the dual frame $\Phi_{\tilde F}$ is $\tilde F=F(F^tF)^{-1}$. Since $\tilde F^t \tilde F=(F^tF)^{-1}$,the eigenvalues of $\tilde F^t\tilde F$ lie in the closed interval $[1/B,1/A]$, and the bounds of $\Phi_F$ are $B^{-1}$ and $A^{-1}$. Since $\tilde F F=I_K$, given the frame expansion coefficients of any vector {\bf z}, the vector can be recontructed using the dual frame operator as follows:
\be
{\bf z}= \tilde F(F{\bf z})=\sum_{k=1}^{N}\langle{\bf z},\varphi_k\rangle\tilde\varphi_k.
\ee

Among the frames, tight frames are the most imporatant since they possess some desirable properties related to reconstruction.The noise sensetivity of a frame depends on the closeness of its bounds. Small perturbations of the frame expansion coefficients can give rise to large reconstruction errors if the frame bounds are far apart.   

\subsection{Filter Banks}

We consider an $N$-channel Filter Bank(FB) with subsampling by the integer factor $M$ in each channel, Perfect Reconstruction(PR) and zero-delay,so that $\hat{x}[n]=x[n]$ where $\hat{x}[n]$ and $x[n]$ denote the reconstructed and the input signal respectively. The transfer functions of the analysis and synthesis filters $H_k(z)$ and $F_k(z)$ ($0\leq k\leq N-1$),with corresponding impulse responses $h_k[n]$ and $f_k[n]$ respectively. In the oversampled case,$(N>M)$, the subband signals are redundant, since they contain more samples(per unit of time)than the input signal. Oversampled FBs offer more design freedom, improved numerical properties and improved noise immunity as compared to critically sampled FBs. The increased design freedom corresponds to the fact that, for a given oversampled analysis FB, there exists a whole class of synthesis FBs providing PR.

The \emph{polyphase decomposition} of the analysis filters $H_k(z)$ reads 
\be
H_k(z)=\sum_{n=0}^{M-1}z^n E_{k,n}(z^M),0\leq k\leq N-1
\ee where
\be
E_{k,n}(z)=\sum_{m=-\infty}^{\infty}h_k[mM-n]z^{-m}
\ee
with $0\leq k\leq N-1,0\leq n\leq M-1$ is the $n$th polyphase component of the $k$th analysis filter $H_k(z)$.The $N\times M$ analysis polyphase matrix ${\mathbf E}(z)$ is defined as $[{\mathbf E}(z)]_{k,n}=E_{k,n(z)}$. The synthesis filters can be decomposed similarly.

If we now consider subband signals $v_k[m]$ corresponding to input signal $x[n]$ and the reconstructed signal $\hat{x}[n]$ and the perturbed subband signals $v'_k[m]=v_k[m]+\Delta v_k[m]$ corresponing to the input signal $x'_k[m]=x_k[m]+\Delta x_k[m]$ and the reconstructed  signal $\hat{x}'_k[m]=\hat{x}_k[m]+\Delta \hat{x}_k[m]$. Using eqn.\ref{frame} and \cite{framefb}, it can be shown that 
\be
A \leq \frac{\parallel \Delta \hat{x}\parallel^2}{\parallel \Delta v\parallel^2}\leq B
\ee
where $A,B$ are the frame bounds for the vectors $\Delta v_k[m]$ over $l^2({Z})$. Thus for a given subband perturbation energy $\parallel \Delta v\parallel^2$, frame bounds $A$ and $B$ provide lower and upper bounds on the resulting reconstruction error energy.The reconstruction error energy is minimized by making $A$ as small as possible and $B \approx A$. Thus it is desirable to have a snug frame. The frame bounds can also be related to the \emph{oversampling factor} $N/M$. Normailzations give \cite{framefb}
\be
A \leq \frac{1}{N/M}\leq B.
\ee
For tight frames, corresponding to a paraunitary FB, we have $ A = B= \frac{1}{N/M}$.

Thus, the energy of the reconstruction error for given subband perturbation energy $\parallel\Delta v\parallel^2$ is here inversely proportional to the oversampling ratio. 

For further details on the relation between frames and FBs, the reader is referred to \cite{framefb}.

\section{Codes over the Real Field ${R}$}       

\subsection{Filter Banks as Convolutional Encoders}
In this chapter we begin with our idea of iterative decoding if codes over the real field for robustness to erasures. Although this idea we believe can be generalised to two dimensions as well, such as for images, only the one dimensional case is dealt with in this report. We use rate $\frac{1}{2}$ codes for encoding the data to begin with. Multirate (rate $1/M$) codes are dealt with later. In earlier works \cite{motwani},the use of oversampled filter banks for robustness against erasures has been studied. Although the connection between filter banks and convolutional encoders is obvious, there seem to exist only a few publications devoted to this sublect, see \cite{xia},\cite{andreas}. The frame bounds are a crucial factor in deciding the reconstruction mean square error when the codewords are quantized. It has been shown \cite{goyal} that if the encoding is done using uniform frames, the reconstruction mean square error is minimized if and only if the frame is tight.

It is seen for DFT codes, in the case of cosecutive erasures, this value can be quite high \cite{dft}. However, if the data is encoded using two encoders as in the case of turbo codes, significant improvements in perforamance were observed in \cite{dft}. It was observed that the eigen spread of the resultant equivalent frame ${\mathbf T_r}$ was smaller. 

\section{The Encoding Procedure - Rate $1/2$ Codes}
Let the input sequence $x(n)$ be blocked into vectors of size $N$ (where $N$ is the size of the interleaver) denoted by ${\mathbf x}$, i.e., ${\mathbf x}= [x(0) x(1) x(2) \dots x(N-1)]^t.$ Let the output sequences of the filter $x_0(n)$ and $x_1(n)$ be combined into one vector ${\mathbf y},$ i.e.,
\be
{\mathbf y}=[x_0(0) x_0(1)x_1(0) x_1(1) x_0(2) \dots x_0(\frac{N-1}{2})x_1(\frac{N-1}{2})]^t
\ee
 and the output sequences of the filter banks $x_0^{\pi}(n)$ and $x_1^{\pi}(n)$ be combined into one vector ${\mathbf y^{\pi}}$,i.e.,
\be
\label{pi}
{\mathbf y^{\pi}}=[x_0^{\pi}(0) x_0^{\pi}(1)x_1^{\pi}(0) x_1^{\pi}(1) x_0^{\pi}(2)\dots x_0^{\pi}(\frac{N-1}{2})x_1^{\pi}(\frac{N-1}{2})]^t.
\ee 
We can then express the output ${\mathbf y}$ of the fiter bank as a linear transform ${\mathbf T}$ acting on the input ${\mathbf x}.$ If the analysis filters are causal and have length $L$, the matrix ${\mathbf T}$ is given by 
\be
{\mathbf T}=\left[\begin{array}{c}
			{\mathbf T_s}\\
			{\mathbf T_{\pi}}\end{array}\right]
\ee
where
\be
\label{big}
{\mathbf T_s}=\left(\begin{array}{ccccccccc}
                    \ddots & & & & & & && \\
                       & h_0(L-1)& h_0(L-2)&h_0(L-3)&\ldots &h_0(0)&0 &0&  \\
		       & h_1(L-1)& h_1(L-2)&h_1(L-3)&\ldots &h_1(0)&0 &0& \\
	        	 & 0 & 0 & h_0(L-1)&\ldots &h_0(2)&h_0(2) &h_0(0)&\\
                        & 0 & 0 & h_1(L-1)&\ldots &h_1(2)&h_1(2) &h_1(0)& \\
                        & & & & & & & & \ddots \end{array}\right)
\ee 
and $\Pi$ is a pseudo random interleaver. ${\mathbf T_{\Pi}}$ is such that its $i$th row is the $j$th row of ${\mathbf T_s}$ where ${\Pi}(i)=j.$ 

In the event of erasures the reconstuction mean square error is given by \cite{motwani}:
\be
MSE_{\sigma^2}= \sigma^2 trace({\bf T}_r^t{\bf T}_r)^{-1}.
\ee
where ${\mathbf T_r}$ is obtained at the receiver (decoder) after some rows from the matrix ${\mathbf T_s}$ and the corresponding rows from ${\mathbf T_{\Pi}}$ have been erased.

\subsection*{Comments}
 As a result of the arbitray nature of the interleaver it is not possible to calcultae the trace of the square matrix $({\bf T}_r^t{\bf T}_r)^{-1}.$ Computer simulations were done for an interleaver size of 150. It was found that $ trace({\bf T}_r^t{\bf T}_r)^{-1}$ is around \emph{five} orders of magnitutde smaller than those obtained using oversampled filter banks in \cite{motwani}. The reconstruction is possible for upto $2L-1$ errors provided they are not in the first $L-1$ rows of ${\mathbf T_s}.$ This is because reconstruction is only possible if the remaining rows of ${\mathbf T_s}$ and ${\mathbf T_{\Pi}}$ form a frame. Physically speaking, there should be enough redundancy in the signal left after erasures for reconstruction to be possible. Even though there were certain \emph{rogue} permutations for which the improvement in trace in not so dramatic they are in the worst cases two orders of magnitude less than the oversampled case. As we know, the performance of Turbo codes is to a large extent governed by the exact permutation at hand. In this case however, the performance was found to be independent of the cycle length of the permutation and hence ${\Pi}.$ It was also found that our procedure of encoding against burst erasures over erasure channels of the type normally encountered in multimedia signalling is better than that in \cite{dft}. 
 
There are certain problems with the simulations though, of which one should be aware. The matrix (\ref{big}) is bi-infinite. The simulations however use finite matrices and filters (in our case of length 9). As a result the orthogonality of the rows of ${\mathbf T_s}$ and ${\mathbf T_{\Pi}}$ is not guaranteed. In fact for a filter of length $L$, the first and last $L+1$ rows are not mutually orthogonal, while the middle $N-2(L-1)$ rows are mutually orthogonal. These create some problems in the reconstruction which we address next.

\section{Reconstruction}
As we are dealing with codes over real fields, the conventional problem of decoding over $GF(2)$ appears as the problem of reconstruction. The channel that we address is the erasure channel and thus our task is to recover the values lost on transmission. The origin of noise in our case would be quantization. This has been a standard problem in Image processing and substantial amount of literature exists in that area. Of them we will choose one particular reconstrucion algorithm for our purposes.  Prior to that however we have to ensure the recoverability of the maximum number of erasures. This is dealt with in the next section.
\section{Recovery}
In this context the best permutation is $\Pi(i)=i+ N/2 (mod N)$, where $N$ is the interleaver size. We can always choose $N$ to be even as it is in our prerogative. This particular problem of decoding over $GF(2)$ appears as the problem of reconstruction. The channel that we address is the erasure channel and thus our task is to recover the values lost on transmission. The origin of noise in our case would be quantization. This has been a standard problem in Image processing and substantial amount of literature exists in that area. Of them we will choose one particular reconstrucion algorithm for our purposes.  Prior to that however we have to ensure the recoverability of the maximum number of erasures. This is dealt with in the next section. In this spirit, we have the following theorem.

{\bf Theorem 1:}The permutation $\Pi(i)=i+ N/2 (mod N)$ can recover upto $N$ consecutive erasures and no more when passed through the system given by Eqn \ref{big}.

Proof: When the interleaver size is N it is evident that the matrix $\mathbf{T}$ is of dimension $2N \times N$. When an information string of length $N$ is operated upon by this $\mathbf{T}$, we obtain the the sequence in Eqn \ref{pi}. With slight abuse of notation we will label this sequence as $\mathbf{y}$ where 
\begin{center}
\be
		 \left(\begin{array}{cc}
                                      y_1 \\
				      y_2 \\
				      \vdots \\ 
                                       y_{N-1}\\       
                                       y_{N} \\
					\vdots\\
					y_{2N -1}\\
					y_{2N}	\end{array}\right)= \left(\begin{array}{ccccccccc}
                       & h_0(L-1)& h_0(L-2)&h_0(L-3)&\ldots &h_0(0)&\ldots &0&  \\
		       & h_1(L-1)& h_1(L-2)&h_1(L-3)&\ldots &h_1(0)&\ldots &0& \\
	        	 & 0 & 0 & h_0(L-1)&\ldots &h_0(2)&h_0(2) &h_0(0)&\\
                        & 0 & 0 & h_1(L-1)&\ldots &h_1(2)&h_1(2) &h_1(0)& \\
			& \vdots & \vdots &  &\ldots & &  & \vdots& \\
			& 0 & \ldots &h_0(L-1)  &h_0(L-2) &\ldots &h_0(0) &\ldots& \\
			& 0 & 0 &\ldots  & \vdots&\ldots &h_0(L-3) &\ldots&\\\end{array} \right)\left(\begin{array}{cc}
                                      x_1 \\
				      x_2 \\
				      \vdots \\ 
                                       x_{N-1}\\       
                                       x_{N} \end{array}\right).
\ee
\end{center}
This will give us a set of $2N$ equations in $N$ variables forming a consistent system. The nature of erasures considered are such that if $y_i$ is lost then so is $y_{i+N}$. However $y_i$ if it contains information about say $y_{\alpha_{i}}$, $\alpha_i$ is a set of $L$ integers from 1,2,\dots,N, then $y_{i+N}$ has information about $y_{\alpha_{i}+N}$. The sum is modulo N. As a result, after the loss of $N$ of these equations, the remaining still form a system of $N$ equations in $N$ variables which can be solved to obtain the values of the coefficients $x_1, x_2,\dots, x_N.$

	Having proved that the permutation $\Pi(i)=i+ N/2 (mod N)$ does ensure recovery of the maximum number of erasures which is possible in this case, {\it i.e.,} $N$ we try to devise a scheme for the rather unlikely scenario of being able to recover after more than $N$ erasures from $2N$ samples. 

\subsection{M - Channel Filter Banks}

	Suppose now that we have to decompose the initial signal into a direct sum of $M$ subspaces. We try to develop the matices analogous to (\ref{big}) on the lines of \cite{cormac}, \S II. B. There will be $M$ filters $h_i$, $ i=1,2,\dots, M-1$. They have to satisfy the basic requirement for projections onto orthonrmal and complete subspaces
\be
\sum_{i=0}^{M-1}{\mathbf H_i^{*}}{\mathbf H_i}={\mathbf I}.
\ee
 In that case the T matrix analogous to \ref{big} will be given by 
\be
\label{mch}
{\mathbf T_s}=\left(\begin{array}{ccccccccc}
                    \ddots & & & & & & && \\
                       	& h_0(L-1)& \ldots&h_0(L-M+1)&\ldots &h_0(0)&0 &0&  \\
		       	& \vdots& \vdots&\vdots&\ldots &\ldots&0 &0& \\
	              	& h_M(L-1)& \ldots&h_M(L-M+1)&\ldots &h_0(0)&0 &0&  \\	 
			& 0 & 0 & h_0(L-1)&h_0(L-2)&h_0(L-3)&\ldots &h_0(0)&\\
                        & 0 & 0 & h_1(L-1)&h_1(L-2)&h_1(L-3)&\ldots &h_0(0)&\\
                        & & & &\vdots & & \vdots& & \ddots \end{array}\right)
\ee
{\it i.e.,} the subsequent sets of $M$ rows are shifted by $M$ each to the right. This system has the following as independent parameters: $L$, $M$ and $N$. To devise the $\Pi$ in the most general setting is a daunting task and is still an open problem. However, we will prove a couple of theorems for some special cases that in essence prove that $M$ channel filter banks may indeed be better choices. 

{\bf Theorem 2:} If the M channel filter bank is such that $L = M$ and $N=rM$ for $r\in {Z}^{+}$, then the system can correct upto a maximum of $2(kM - 1)$ errors, such that $(k-1)M<N/2<kM$. 

Proof: From the matrix \ref{mch} and applying the permutation $\Pi(i)=i+ N/2 (mod N)$ on it, we get $2N$ equations such that\ben
y_1,\dots,y_M &\}& x_1,\dots,x_M \\ \nonumber
y_{M+1},\dots,y_{2M} &\}& x_{M+1},\dots,x_M \\ \nonumber
\vdots \\ \nonumber
y_{(r-1)M+1},\dots,y_{rM}&\}&x_{(r-1)M+1},\dots,x_{rM}\\ \nonumber
y_{N+1},\dots,y_{N+M} &\}& x_{N/2 + 1},\dots, x_{N/2+M}\\\nonumber
\vdots \\ \nonumber
y_{2N-M+1},\dots,y_{2N} &\}&x_{N/2 -M+1},\dots,x_{N/2}.
\een 
W.l.o.g, we can count how many erasures one can correct beginning from the first one. The channel is still such that if we lose $y_i$ then we loose $y_{i+N}$. Find the $k$ such that $(k-1)M<N/2<kM$, then it is evident that we can recover erasures upto $(k-1)M$. To find if we can recover after more erasures, we see that we cannot recover if we loose upto $kM$ as in that case there will be no equation containing at least $x_1$. This can be seen as follows. When we loose $y_1,\dots,y_{kM}$, the loss is $x_1,\dots,x_{kM}$. Correspondingly we also loose $y_{N+1},\dots,y_{kM+N}$ whence the loss is $x_{N/2+1},\dots,x_{N/2 + kM}$. Since all sums are modulo $N$ and $kM > N/2$, $(N/2 + kM) \geq 1$. 
Hence we should at least leave one equation from the last set of $M$ equations. Hence the number is $kM -1$ out of $N$. Thus it can correct upto $2(kM-1)$ erasures in $2N$ symbols. 

We will be interested in the effect of the amount of channeling and its effect on recoverability in the above scheme. We have the following result in that context. 

{\bf Corollary 1:} With the above scheme, one can recover from more than half losses.
	
If $N=150$ which is mapped to 300 points and for $M=4=L$, we can recover after a loss of 150 elements. For  $M=8=L$, we can recover after a loss of 158 elements.For $M=16=L$, we can recover after a loss of 158 elements. For $M=32=L$, we can recover after a loss of 190 elements.

Of course, this corroborates the fact that the longer the filter the better but since we assumed $L=M$, the structure will become unwieldy. Also computational complexity will be a challenge.

In general $L>M$ for all practical purposes. In this context we have another partial result.

{\bf Theorem 3:} If the M channel filter bank is such that $L = 2M$ and $N=rM$ for $r\in {Z}^{+}$, then the system can correct upto a maximum of $N+2M$ errors.

Proof: Similar to the proof of the last theorem we have 
\ben
y_1,\dots,y_M &\}& x_1,\dots,x_{2M} \\ \nonumber
y_{M+1},\dots,y_{2M} &\}& x_{M+1},\dots,x_{3M} \\ \nonumber
\vdots \\ \nonumber
y_{(r-1)M+1},\dots,y_{rM}&\}&x_{(r-1)M+1},\dots,x_{rM+M}\\ \nonumber
y_{N+1},\dots,y_{N+M} &\}& x_{N/2 + 1},\dots, x_{N/2+2M}\\\nonumber
\vdots \\ \nonumber
y_{2N-M+1},\dots,y_{2N} &\}&x_{N/2 -M+1},\dots,x_{N/2+M}.
\een

	However we note that there are overlaps in the sequence. Any set of $M$ $x$'s appear twice in consecutive groups of $y$'s each of size $M$. Hence using the argument of the previous Theorem we see that we can recover from  more than $kM$ erasures where $(k-1)M<N/2<kM$. Let $\Delta = N/2 - (k-1)M$, then we can loose upto $kM + \Delta = M+ N/2$ elements and still recover the initial sequence. 

Thus from a sequence of $2N$ elements one can recover if there is an erasure of upto $N + 2M$ elements.

\subsection*{Comments}

	The results of the last two sections prove the existence and recoverability of the initial information sequences fora variety of systems and provide the infrastructure to obtian results for a multitude of other systems. However the existence of recoverability does not the solve the engineering problem of decoding. What is needed is a synthesis Filter Bank that can recover the initial information string given any arbitrary pattern of erasures subject to the bounds provided by the above theorems. This is a formidable task and appears tough to solve in the present scenario. 
	We however look at a simple algorithm that may  provide the motivation for a reconstruction procedure.

\subsection{Alternating Projection Theorem}
The method of generalized image restoration by Alternating Orthogonal Projections was pioneered by Youla \cite{youla}. The concept of the method is very simple in mathematical terms. Given a Hilbert space ${\mathcal H}$ with  elements $f,g,h$ etc, a zero vector $\phi$ and an inner product, we consider an element $f\in {\mathcal H}$ belonging to a \emph{known} closed linear manifold(CLM) ${\mathcal P}_b$. We are only given its projection $g={\mathcal P}_a f$ onto the known CLM  ${\mathcal P}_a$. The algorithm to reconstruct $f$ from $g$ is given by \cite{youla} as 
\be
f_{k+1} = g + {\mathcal Q}_a{\mathcal P}_b f_k
\ee
where ${\mathcal Q}_a$ is the orthogonal complement of ${\mathcal Q}$ and ${\mathcal P}_a$ is the orthogonal complement of ${\mathcal P}.$  
The method always works in the absence of noise and in the presence of noise it works under certain restrictions on the nature of the noise and the from of ${\mathcal P}$ and ${\mathcal Q}$. 
 
In our case we can use a projection like procedure. But as it turns out in the absence of any quantization noise, we do not need any recursion. Our algorithm get back the erased values in one run only. The manifolds in our case are as follows:\\
${\mathcal P}_a$: space spanned by the remaining (not lost) rows of ${\mathbf T_s},$ \\ 
${\mathcal P}_b$: space spanned by the remaining (not lost) rows of ${\mathbf T_{\Pi}},$ \\
${\mathcal Q}_a = {\mathcal P}_a^{\perp}$: space spanned by the lost rows of ${\mathbf T_s},$ \\
${\mathcal Q}_a = {\mathcal P}_a^{\perp}$: space spanned by the lost rows of ${\mathbf T_{\Pi}}.$ \\
 We also note that either 
\be
 {\mathcal P}_a^{\perp}\subset {\mathcal P}_b 
\ee 
 or 
\be
{\mathcal P}_b^{\perp}\subset {\mathcal P}_a
\ee
\label{second}
for any reconstruction to be possible. 
Then our recontruction relation is 
\be
x = {\mathcal P}_a x + {\mathcal P}_a^{\perp}{\mathcal P}_b x.
\ee
\label{reco}
where $x$ is the information vector in Fig:4.1.
In case Eqn: (5.3) holds the relation is
\be
x = {\mathcal P}_b x + {\mathcal P}_b^{\perp}{\mathcal P}_a x.
\ee 
We actually work with the cofficients of the vector $x$ rather than the vector itself for the sake of computational efficiency. Thus all we need to obtain on the RHS of Eqn: (5.4) are the coefficients of ${\mathcal P}_b^{\perp}{\mathcal P}_b x.$ The $j$th coefficient of $ {\mathcal P}_b^{\perp}{\mathcal P}_b x$ can be given by:
\ben
[{\mathcal P}_a^{\perp}{\mathcal P}_b x]_j = \langle b_j, \sum_i\langle x,e_i\rangle e_i \rangle
					   = \sum_i \langle x,e_i\rangle \langle b_j,e_i\rangle.
\een
where $b_j$ span ${\mathcal P}_a^{\perp}$ and $e_i$ span ${\mathcal P}_b.$

Once again due to the nature of the matrices and the spaces ${\mathcal P}_a$,${\mathcal P}_b$ and their orthogonal complements, the reconstruction procedure cannot be verified analytically. Numerical simulations agree with our predictions but there are certain problems. The bases $b_j$ and $e_i$ are not complete due to the finiteness of the matrix (\ref{big}). Hence we need to do zero-tailing and then look at the reconstruction.

\section{Conclusion and Open Problems}

	Since the result of Shannon in 1948, the continued aim of mathematicians and communication engineers has been to design codes that approach the channel capacity. Block codes existed before Shannon and inspite of their beautiful mathematical structure, they failed to excite communication engineers except in a very few cases, like the Reed Solomon Codes. The emergence of Convolutional codes was more of an engineering solution to the challenge of attaining the capacity. Forney provided an algebraic structure and these codes went closer to the capacity limit given by Shannon. It was not until 1993 however when  C. Berrou, A. Glavieux,P.Thitimajshima touched upon the Shannon limit asymptotically. Out of a computer search fell the Turbo codes. They were found to achieve the limit. However, lack of analytical results in this context dissapoint a lot of people with a mathematical bent of mind. It was evidently a dream come true for communication engineers.
	Berrou, Glavieu, Thitimajshima were quite right in developing turbo codes over $GF(2)$. It was a natural question to ask if these codes can be extended to continuous fields. However, until upto 10 years of  the emergence of Turbo codes this question was not asked, lest answered. 

	This has been the motivation behind this article. We have tried, in a small way., to develop the idea of Turbo Codes over ${R}$. This, we believe, will be a very active area of reasearch involving the concepts of coding theory, frame theory and signal processing. Each of these are huge disciplines by themselves.

	We have used techniques from these fields and some very basic mathematics to arrive at our results. Simulations have provided initial corroboration to our ideas. We have come up with some basic results proved in a limited setup. However these results can be used as stepping stones to obtain much deeper results about Turbo Codes over the real field.

	These ideas can also be extended to higher dimensional systems, in the least to the case of images, which is ${R}^2$. We end with some of the \emph{open problems} in this area :

\begin{itemize}
\item{To obtain the synthesis FB in 2 channel case for an arbitrary erasure pattern.}
\item{The behaviour of the code in the presence of Channel noise other than quantization noise.}
\item{To obtain general expression for maximum recoverability in M channel case.}
\item{To obtain the synthesis FB in M channel case for an arbitrary erasure pattern.}
\item{To explore the possibility of this scheme in the 2D case and obtain corresponding results as above.}
\end{itemize}

\section*{Acknowledgements}
        I would like to take this opportunity to express my sincere and earnest regards for Dr. Ravi Motwani.

\end{document}